\documentclass[article]{amsart}
\usepackage{url}
\usepackage{mathrsfs}

\def\1{1\!\!\!1}

\newcommand{\ceil}[1]{\left\lceil #1 \right\rceil}

\newcommand{\md}[1]{(\text{mod}\ #1)}

\newcommand{\Acal}{\mathcal{A}}
\newcommand{\Bcal}{\mathcal{B}}
\newcommand{\Ccal}{\mathcal{C}}
\newcommand{\sg}{\sigma}

\newcommand{\eps}{\epsilon}
\newcommand{\dl}{\delta}

\newcommand{\mb}{\mathbb}

\newcommand{\lf}{\left}
\newcommand{\rt}{\right}
\newcommand{\lm}{\limits}
\newcommand{\st}{\substack}

\newcommand{\lb}{\lambda}

\newcommand{\fr}{\frac}

\newcommand{\vp}{\varphi}

\newcommand{\mc}{\mathcal}

\newcommand{\Z}[1]{\left( \mb{Z}/#1 \mb{Z} \right)^\times}

\newtheorem{thm}{Theorem}[section]

\newtheorem{lem}[thm]{Lemma}

\newtheorem*{rem}{Remark}

\begin{document}
\title{Product of three primes in large arithmetic progressions}
\author{R. Balasubramanian}
\address{
Institute of Mathematical Sciences, Taramani Chennai India-600113
\\and Homi Bhabha National Institute, Training School Complex\\
Anushakti Nagar, Mumbai, India-400094.\\
}
\email{balu@imsc.res.in}

\author{O. Ramar\'e}
\address{
CNRS / Institut de Math\'ematiques de Marseille\\
Aix Marseille Universit\'e, U.M.R. 7373\\
Site Sud, Campus de Luminy, Case 907\\
13288 MARSEILLE Cedex 9, France\\
}
\email{olivier.ramare@univ-amu.fr}

\author{P. Srivastav}
\address{
Institute of Mathematical Sciences, Taramani Chennai India-600113
\\and Homi Bhabha National Institute, Training School Complex\\
Anushakti Nagar, Mumbai, India-400094.\\
}
\email{priyamvads@imsc.res.in}



\begin{abstract}
  For any $\epsilon>0$, there exists
  $q_0(\epsilon)$ such for any $q\ge q_0(\epsilon)$ and any invertible residue
  class $a$ modulo~$q$, there exists a natural number
  that is congruent to $a$ modulo~$q$ and that
  is the product of exactly three primes, all of which are below
   $q^{\frac{3}{2}+\epsilon}$.
  If we restrict our attention to odd moduli $q$ that do not have prime
  factors congruent to~1 mod~4, we can find
  such primes below $q^{\frac{11}{8}+\epsilon}$.
  If we further
  restrict our set of moduli to prime $q$ that are such that
  $(q-1,4\cdot7\cdot11\cdot17\cdot23\cdot29)=2$, we can find
  such primes below $q^{\frac{6}{5}+\epsilon}$. Finally, for any $\epsilon>0$, there exists
  $q_0(\epsilon)$ such that when $q\ge q_0(\epsilon)$, there exists a natural number
  that is congruent to $a$ modulo~$q$ and that
  is the product of exactly four primes, all of which are below
  $q(\log q)^6$.
\end{abstract}
\maketitle{}

\section{Introduction and results}

In this paper we investigate the representation of reduced residue
classes modulo~$q$ by a product of exactly three small primes, and
exceptionally by four of them in Theorem~\ref{mainthmasympter}.
We develop the approach initiated
in~\cite{Ramare-Walker*16}, and followed up
in~\cite{Ramare-Serra-Srivastav*18}, borrowing several results from
the latter paper.

We input three new arguments: some vertical Brun-Titchmarsh
inequalities (e.g. Theorems~\ref{iwa} and \ref{mika}), the usage of
P${}_2$-numbers, see Theorem~\ref{P2} and some additive combinatorics
of sum-free sets in Section~\ref{Sac}.

When using
Kneser's theorem, we encounter two hurdles: the
need to prove that the set of classes of small primes modulo~$q$ is
dense enough, this is the job of the different versions of the
Brun-Titchmarsh inequalities we employ, and the  need to show that
these primes do not stay in some union of cosets of some subgroup of
small index in~$\Z{q}$. By following the strategy of P.D.T.A. Elliott
in~\cite{Elliott*71a}, we prove that subgroups of index~5 contain
small primes (we prove more, see Theorem~\ref{T1}). Finally
Theorem~\ref{P2} asserts that every coset of a small index
subgroup contains either a small prime or a product of two of
them. This is a consequence of the weighted sieve when using the
approach we already put to work for the coset Brun-Titchmarsh inequality
in~\cite{Ramare-Serra-Srivastav*18}. The surprise is that, though we
seem to be using the same kind of sieve argument as when bounding the
density from above, the additive
consequences are distinct. The additive combinatorial problem that
emerges is investigated in Section~\ref{Sac}.
It relies on the combinatorics of sum-free sets.

\begin{thm}
  \label{mainthmasymp}
  Let $\epsilon>0$. There exists $q_0$ such that
  for all $q\ge q_0$ and for all inver\-ti\-ble residue classes $a$ modulo~$q$,
  there exists an integer congruent to $a$
modulo~$q$ that is the product of exactly three primes, all of
which are below~$q^{\frac{3}{2}+\epsilon}$. When $q$ is cube-free, we
can find such primes below $q^{\frac{4}{3}+\epsilon}$.
\end{thm}
The main interest in the above result is that $\frac{3}{2}<2$ while
it is unknown whether there is such a small prime in any given
arithmetic progression modulo~$q$, the best bound for Linnik's
constant under the Generalized Riemann Hypothesis being $2+\epsilon$.

We encounter two obstructions
during the proof: the parity phenomenom and the large subgroups
problem (arising from the usage of Kneser's theorem). We can avoid the
second one by specializing the modulus $q$ to a well-behaved family.
\begin{thm}
  \label{mainthmasympquatro}
  Let $\epsilon>0$. There exists $q_0$ such that
  for every modulus $q\ge q_0$ all whose prime factors are congruent
  to $3$ modulo~4, and for all 
  inver\-ti\-ble residue classes $a$ modulo~$q$, there exists an
  integer congruent to $a$ modulo~$q$ that is the product
  of exactly three primes, all of which are
  below~$q^{\frac{11}{8}+\epsilon}$.
\end{thm}
  The same applies to integers $q$ of the form $q=4q'$ or $q=8q'$
  when all the prime factors of $q'$ are congruent
  to $3$ modulo~4.
  Notice for comparison that $11/8=1.375$.
\begin{thm}
  \label{mainthmasympbis}
  Let $\epsilon>0$. There exists $p_0$ such that
  for every prime $p\ge p_0$ such that
  $(p-1,4\cdot7\cdot11\cdot17\cdot23\cdot29)=2$ and for all 
  inver\-ti\-ble residue classes $a$ modulo~$p$, there exists an
  integer congruent to $a$ modulo~$p$ that is the product
  of exactly three primes, all of which are
  below~$p^{\frac{6}{5}+\epsilon}$.
\end{thm}
We can reach a smaller exponent for a general modulus by taking
products of \emph{four} primes, rather than of \emph{three}.
\begin{thm}
  \label{mainthmasympter}
  There exist $q_0$ such that for all
  modulus $q\ge q_0$ and for all
  inver\-ti\-ble residue classes $a$ modulo~$q$, there exists an
  integer congruent to $a$ modulo~$p$ that is the product
  of exactly four primes, all of which are
  below~$q(\log q)^6$.
\end{thm}
In \cite{Erdos-Odlyzko-Sarkozy*87}, P.~Erd\H{o}s, A.~Odlyzko and
A.~S\'ark\"ozy proved that under the Generalized Riemann Hypothesis,
we can find such a product of only two primes, and that under a weaker
but still unproven hypothesis, we can find such a product of~three
primes.  It is unclear to us how their method would work for products
of four primes in case of the presence of a Siegel zero.

\bigskip
\noindent
{\bf Thanks} are due to the referee for his/her careful
reading of this paper.

\section{Primes in large subgroups modulo~$q$}
\label{Elliott}


\begin{lem} 
\label{L1chi} \phantom{mm}
\begin{itemize}
    \item Given $\eps>0$, there exists $c(\eps)>0$ such that, for
  all real characters $\chi$ modulo~$q$, we have
    $|L(1, \chi)| > c(\eps) \, q^{-\eps} $.
    \item There exists a constant $c>0$, such that for all complex
      characters $\chi$ modulo $q$, we have $ |L(1, \chi)| \ge c/\log q$. 
 \end{itemize}
\end{lem}

\begin{lem}[Burgess]
\label{Burgess}
Let $\chi$ be a non-trivial character to modulus $q$. Then for all $X$,
$H\ge1$, and for all $r\in\{1,2,3\}$, we have the following: Given
$\eps>0$, there exists $c(\eps)$, such that
\begin{equation*}
 \sum\lm_{X < n \leq X + H} \chi(n)  \le c(\eps)  H^{1-\fr{1}{r}} q^{\fr{r+1}{4 r^2} + \eps}.
\end{equation*}
Further, if $q$ is cube-free, one can take $r$ to be any natural number.
\end{lem}

The following is a consequence of the Burgess bounds.
\begin{lem}For every character modulo~$q$, we have
\label{chiB}
\begin{itemize}
    \item If $y \geq q^{1/3 + \eps}$, there is a $\dl>0$ such that
    $$ \sum\lm_{n \leq y} \fr{\chi(n)}{n} = L(1,\chi)  + O(q^{-\delta}) .$$ 
    The factor $1/3$ can be replaced by $1/4$ when $q$ is cube-free.
    
    \item Let $r \in \{ 1,2,3 \}$ and suppose that $0<\sg < 1 - 1/r$. Then
    $$ \sum\lm_{n \leq y} \fr{\chi(n)}{n^{s}} \ll \fr{y^{1 - \sg - \fr{1}{r} }}{1 - \sg - 1/r} \ q^{\fr{r+1}{4r^2} + \eps} \, |s|  . $$
    Again, we can take any $r \geq 1$ when $q$ is cubefree. 
\end{itemize}
\end{lem}

\begin{thm}
\label{T1}
Let $q$ be sufficiently large and
$H \subseteq G_q = \lf( \mb{Z}/q \mb{Z} \rt)^*$ be a subgroup of index
$Y$. Then, there exists a prime
$p \leq c_0(\eps, Y) \cdot q_0^{\fr{Y-1}{3} + \eps}$, such that
$p \md{q} \in H $. Further, if $q$ is cube-free, one can replace the
exponent $(Y-1)/3$ by $(Y-1)/4$.
\end{thm}

In
\cite{Elliott*71a}, P.D.T.A. Elliott proved the same result when $H$
is the set of invertible $Y$-th powers modulo~prime $q$. Our
proof follows along the similar lines.  When $Y=2$, the paper
\cite{Vinogradov-Linnik*66} of A.I. Vinogradov \& Y. Linnik tells us
that there is a prime $\ll q^{1/4+\varepsilon}$ in a quadratic
subgroup modulo~$q$ \emph{when $q$ is prime}, a condition that has
been removed by P.~Pollack in~\cite{Pollack*14}.
In the case $Y\in\{3,4\}$ and $q$ is prime, pertinent bounds may be
found in \cite{Nagell*52} by T. Nagell and \cite{Kubra-Pollack*19} by
  B. K\"{u}bra and P. Pollack.

\begin{rem}
To be more precise, the bound we obtain in Theorem \ref{T1} is of the form $c_1(\eps) \cdot \lf(Y^3 q_0^{1/3 + \eps}\rt)^{Y-1}$.
\end{rem}

To unfold the proof, we need to introduce the Dirichlet series given
in~\eqref{Fdef}, and this requires some preliminaries.
 The subgroup $\{\chi \in \hat{G}_q: \chi
\vert_H = 1\}$
may be identified with characters on $G_q/H$. We further set $K = G_q/H$.
Let $x > 1$ be such that none of the primes $p \leq x$ falls
into the
subgroup $H$ modulo $q$,

We define a multiplicative function $g(n)$as follows. Let
\begin{equation}
\label{gdef}
g(p) = \fr{1}{Y} \sum\lm_{\chi \in \hat{K}} \chi(p) =
\begin{cases} 1, & p \in H, \\ 0, & \text{otherwise,}  \end{cases}
\end{equation}
and $g(p^\alpha)=0$ for all primes $p$ and $\alpha\ge2$.
Our hypothesis on $x$  implies that 
\begin{equation}
  \label{ghyp}
  g(p) = 0, \quad \text{for all primes } \ p \leq x .
\end{equation}
We let $g(1)=1$. The Dirichlet series of $g(n)$ is given by
\begin{equation}
  \label{fdef}
  f(s) =
  \sum\lm_{n = 1}^{\infty} \fr{g(n)}{n^s}
  = \prod\lm_p \lf( 1 +
  \fr{g(p)}{p^s} \rt)
  = \prod\lm_p \lf( 1 + \fr{1}{Y}\sum\lm_{\chi \in \hat{K}} \frac{\chi(p)}{p^s} \rt) .
\end{equation}
We find that 
\begin{equation*}
f(s)^Y = \omega(s) \prod\lm_{\chi \in \hat{K}} L(s,\chi),  
\end{equation*}
where $\omega(s)$ is Dirichlet series that is analytic for $\sg>1/2$. In fact
\begin{equation}
\label{omega}
\omega(s) = \prod\lm_p  \lf[ \lf(1 + \fr{g(p)}{p^s} \rt)^Y \prod\lm_{\chi \in \hat{K}}\lf( 1 - \fr{\chi(p)}{p^s} \rt) \rt] .
\end{equation}

Let 
\begin{equation}
\label{gtilde}
\tilde{g} = \underbrace{g * \dots * g}_{Y \text{ times}}, \quad \text{so that} \ \, f(s)^Y = \sum_n \fr{\tilde{g}(n)}{n^s} .
\end{equation}
From \eqref{ghyp}, it follows that $\tilde{g}(p) = 0$ for all primes $p \leq x$ and therefore $\tilde{g}(n) = 0$ for all $1<n \leq x$, except for $\tilde{g}(1)=1$. 

\medskip

We now aim to obtain an asymptotic formula for the partial sums
related to $\tilde{g}(n)$ using the Perron's formula. Since we are
only interested in the sum of $\tilde{g}(n)$ with $n$ not exceeding
$x$, we find it convenient to consider the Dirichlet series obtained
by truncating the L-functions at the parameter $x$ (rather than
considering $f(s)^Y$)

\begin{equation}
\label{Fdef}
F(s) = \omega(s) L(s,\chi_{0,q})(s) \prod\lm_{\st{\chi \in \hat{K} \\ \chi \neq \chi_0}} 
\lf( \sum\lm_{n \leq x} \fr{\chi(n)}{n^s} \rt).
\end{equation}

With this notation in hand, we may unfold our series of lemmas.

\begin{lem}
\label{omegaB}
For $1/2< \sg \leq 1$, we have 
$
\zeta(2\sg)^{-2Y} \leq |\omega(s)| \leq \zeta(2\sg)^{2 Y}
$.
\end{lem}

\begin{proof}
We have
\begin{align*}
\log \omega(s) 
&= Y \sum\lm_p \log \lf( 1 + \fr{g(p)}{p^s} \rt)  +   \sum\lm_{\chi \in \hat{K}} \sum\lm_p \log \lf( 1 - \fr{\chi(p)}{p^s} \rt) \\
&= Y \sum\lm_{j=1}^{\infty} \fr{(-1)^{j-1}}{j} \sum\lm_p \fr{g(p)^j}{p^{js}} - \sum\lm_{j = 1}^{\infty} \fr{1}{j} \sum\lm_p \fr{1}{p^{js}} \sum\lm_{\chi \in \hat{K}}  \chi(p)^j.
\end{align*}
The contribution from $j = 1$ cancels out. Now, use the fact that $g(p)$ and $\chi(p)$ are bounded by $1$, to obtain
\begin{align*}
| \log \omega(s)| 
&\leq 2 Y \sum\lm_{j = 2}^{\infty} \fr{1}{j} 
\sum\lm_p \fr{1}{p^{j \sg}} 
\leq 2Y \sum\lm_{j=1}^{\infty} \lf( \fr{1}{2j} \sum\lm_p \fr{1}{p^{2\sg j}} + \fr{1}{2j + 1} \sum\lm_p \fr{1}{p^{\sg(2j + 1)}} \rt) \\
&\leq 2Y \log \zeta(2\sg).
\end{align*}
\end{proof}

\begin{lem}
  \label{FBound}
  Let $s=\sigma+it$. When $1/2<\sigma=\Re s<2/3$, we have
  \begin{equation*}
    F(s) 
    \ll \lf( \zeta^2(2\sg) \fr{x^{2/3 - \sg}}{2/3 - \sg} \, q_0^{1/9 +
      \eps} \rt)^{Y-1}
    (1+|t|)^{Y-\fr{1}{2}} . 
  \end{equation*}
\end{lem}

\begin{proof}
Using Lemma \ref{chiB}(b) with $r = 3$ and the bound $$|L(s,\chi_{0,q})(s)| \ll q_0^{\eps} |t|^{\fr{1 - \sg}{2}} \ll q_0^{\eps} |t|^{1/2}, $$ 
as well as the bound for $\omega(s)$ from Lemma \ref{omegaB}, we find that
\begin{align*}
F(s) 
\ll \lf( \zeta^2(2\sg) \fr{x^{2/3 - \sg}}{2/3 - \sg} \, q_0^{1/9 + \eps} \rt)^{Y-1} (1+|t|)^{Y-\fr{1}{2}} . 
\end{align*}
The proof of the lemma is complete.
\end{proof}

\begin{proof}[Proof of Theorem \ref{T1}]
Let $\tilde{g} = \underbrace{g * \dots * g}_{Y \text{ times}}$ be as in \eqref{gtilde}. We have seen that $\tilde{g}(1) = 1$ and $\tilde{g}(n) = 0$ for all $1 < n \leq x$ owing to the hypothesis \eqref{ghyp}. Therefore
\begin{equation}
\label{tg1}
S = \fr{1}{\ell!}\sum\lm_{n \leq x} \tilde{g}(n) \lf( \log \fr{x}{n} \rt)^\ell \ = \fr{\lf( \log x \rt)^\ell}{\ell!},
\end{equation}
for any integer $\ell \geq 0$. 

On the other hand, we use the following version of Perron's formula (assume $y \neq 1$ and $c>1$):
$$
\fr{1}{2 \pi i} \int_{c - i \infty}^{c + i \infty} \fr{y^s}{s^{\ell + 1}} \, ds
= \begin{cases} \dfrac{(\log y)^\ell}{\ell!}, & y > 1,\\ 0, & \text{otherwise}.  \end{cases}
$$
As a result and with $\ell=Y+1$, we find that
\begin{equation*}
S = \fr{1}{2 \pi i} \int_{2 - i \infty}^{2 + i \infty} F(s) \fr{x^s}{s^{\ell + 1}} \, ds .
\end{equation*}
Moving the line of integration to $\Re(s) = \sg$ ($1/2 < \sg < 2/3$) and collecting the residue at $s = 1$, we get
\begin{equation}
\label{Sfinal}
S = \text{Res}_{s=1}F(s) \fr{x^s}{s^{\ell + 1}} 
+  \frac{1}{2i\pi}\int\lm_{\Re(s) = \sg} F(s) \fr{x^s}{s^{\ell + 1}} \, ds.
\end{equation}
The main term above is controlled by Lemma~\ref{L1chi}:
\begin{align*}
\text{Res}_{s=1} F(s) \fr{x^s}{s^{\ell + 1}} 
  &=  x \omega(1) \fr{\vp(q)}{q}
    \prod\lm_{\st{\chi \in \hat{K} \\ \chi \neq \chi_{0,q}}} \lf( \sum\lm_{n \leq x} \fr{\chi(n)}{n} \rt)\\
&\gg  x \fr{\zeta(2)^{-2Y}}{\log \log q} \prod\lm_{\st{\chi \in \hat{K} \\ \chi \neq \chi_{0q}}} \lf( L(1, \chi) + O(q^{-\dl})  \rt) \gg x q^{-\eps},
\end{align*}
where we have used $\fr{\vp(q)}{q} \gg \fr{1}{\log \log q}$ and the
lower bound for $L(1, \chi)$ from Lemma \ref{L1chi}. The implied
constant depends on $Y$ and $\eps$.  Choose
\begin{equation}
\label{sigmaC}
\sg = 2/3 - \dl,  \quad \text{with} \ \, \dl = \fr{1}{\log q}.
\end{equation}
Using \eqref{FBound}, the second term in \eqref{Sfinal} is at most
\begin{align*}
&\ll x^{2/3 - \dl} \lf(\zeta^2(4/3 - 2\dl) \, \dl^{-1} x^{\dl} q^{1/9+\eps} \rt)^{Y-1}  \int\lm_{1}^{\infty} \fr{dt}{|t|^{3/2}} \\
&\ll x^{2/3 + \dl (Y-2)} \lf(13 \, \dl^{-1} q^{1/9+\eps} \rt)^{Y-1}
                                                                                                                                     \ll x^{2/3 + \fr{Y}{\log q_0}} \lf(13 q^{1/9 + \eps} \log q  \rt)^{Y-1}
                                                                                                                                     \ll x q^{-2 \eps},
\end{align*}
as soon as $x \gg q^{\fr{Y-1}{3} + \eps}$. Hence
\begin{equation*}
S \gg  x q^{-\eps}.    
\end{equation*}
On the other hand, we have from \eqref{tg1}, that 
$S = \fr{(\log x)^{\ell}}{\ell!} = \fr{(\log x)^{Y+1}}{(Y + 1)!} \ll x^{\eps}$,
whenever $Y = o(\log x)$ (which is true since $\fr{Y}{\log x} \ll \fr{1}{\log q}$). This leads to a contradiction and completes the proof. 
\end{proof}

\section{Almost primes in cosets of large subgroups modulo~$q$}

 \begin{lem}[Heath-Brown \cite{Heath-Brown*80a}, Petrow \& Young \cite{Petrow-Young*20}]
   \label{HB}
   Let $\chi$ be a Dirichlet character of conductor $r>1$. When
   $1/2\le \sigma=\Re s$, we have
   \begin{equation*}
     L(s,\chi)\ll_\varepsilon (r(|t|+1))^{\frac{3(1-\sigma)^+}{8}+\varepsilon}
   \end{equation*}
   valid for any $\varepsilon>0$. Here
   $(1-\sigma)^+=\max(0,1-\sigma)$. When
   $r$ is cube-free, one may replace $3/8$ by $1/3$.
 \end{lem}

\begin{thm}
  \label{P2}
  There is a constant $C>0$ such that,
  for every subgroup $H$ of $\Z{q}$ of index $Y$ and every coset $b\cdot
  H$, there exists a $P_2$-number 
  that is of size not more than $C\cdot Y^{9/2}q^{0.768}$, where $C$ is
  some effective constant. When $q$ is cube-free, this size may be
  reduced to $C\cdot Y^{9/2}q^{0.683}$.
\end{thm}
By a ``$P_2$-number'', we mean an integer that has
at most two prime factors.
In \cite{Heath-Brown*78}, R. Heath-Brown proves that, for any $a$ prime to
$q$, there is a
$P_2$-number of size $\ll q^{2-0.035}$ congruent to $a$ modulo~$q$,
where the implied constant does not depend on~$a$. We obtain an
exponent $<1$ for the cosets.
\begin{proof}
  Let us use again the setting of the previous section. We start with
  a subgroup $H\subset G_q=\Z{q}$ of (small) index $Y$. Consider
  $\{\chi \in \hat{G}_q: \chi \vert_H = 1\}$.  This can be identified
  with characters on $K=G_q/H$.  Let $b\in G_q$.  To find a small
  $P_2$-number in $b\cdot H$, we shall be using a special case of the
  main result of G.~Greaves in \cite{Greaves*86} (or Chapter~5 of
  \cite{Greaves*01}) which we now state. Let $\Acal\subset [1,X^*]$ be
  a subset of integers. Assume that
  \begin{equation}
    \label{defR}
    \#\{a\in\Acal, \text{$d$ divides
      $a$}\}=\frac{X\rho(d)}{d}+R(\Acal,d)\quad
    \text{with}\ \rho(d)
    =
    \begin{cases}
      1&\text{when $(d,q)=1$},\\
      0&\text{otherwise.}
    \end{cases}
  \end{equation}
  Let $y$ be a \emph{level of distribution}, i.e. a real number such
  that
  \begin{equation}
    \label{eq:2}
    \sum_{d\le y}3^{\omega(d)}|R(\Acal,d)|\ll \frac{X}{\log^2X}.
  \end{equation}
  Then, if $X^*\le y^g$ where $g=2-0.004456$, then there exists
  $a\in\Acal$ having at most two prime factors.

  Define $X^*=\frac{q}{\varphi(q)}YX$ and
  \begin{equation}
    \label{eq:3}
    \Acal=\bigl\{a: a\le X^*, \text{$a$ falls inside $b\cdot H$}\bigr\}.
  \end{equation}
  Let us verify the conditions of the result of Greaves.
 
  We decompose the
  characteristic function of $b\cdot H$ through multiplicative
  characters, i.e. we write
\begin{equation*}
  \1_{b\cdot H}=\frac{1}{Y}\sum_{\chi\in \hat{K}}\overline{\chi}(b)\chi.
\end{equation*}
On defining 
$A_d(s)=\sum_{d|n\ge1}\1_\Acal(n)/n^s$, the truncated Perron summation
formula gives us
\begin{equation*}
  \sum_{\substack{n\in\Acal,\\
      d|n}}1=\frac{1}{2i\pi}\int_{\kappa-iT}^{\kappa+iT}A_d(s)\frac{X^{*s}ds}{s}+O\biggl(
  \frac{X^*\log^2X^*}{dT}\biggr)
\end{equation*}
provided that $1\le T^2\le X^*/d$ and where $\kappa=1+1/(\log X^*)$. This
expands into
\begin{equation*}
  \sum_{\substack{n\in\Acal,\\
      d|n}}1=\frac{1}{Y}\sum_{\chi\in \hat{K}}\overline{\chi}(b) \chi(d) \frac{1}{2i\pi}\int_{\kappa-iT}^{\kappa+iT}L(s,\chi)\frac{X^{*s}ds}{d^ss}+O\biggl(
  \frac{X^*\log^2X^*}{dT}\biggr).
\end{equation*}
We recall Lemma~\ref{HB}, and deduce that
\begin{equation*}
  \sum_{\substack{n\in\Acal,\\
      d|n}}1-\frac{1}{Y}\frac{\varphi(q)}{q}\frac{X^*}{d}\1_{(d,q)=1}
  \ll
  \frac{X^*\log^2X^*}{dT}+\frac{(qT)^{3/16+\varepsilon}\sqrt{X^*}}{\sqrt{d}}.
\end{equation*}
 With $R(\Acal,d)$ defined in~\eqref{defR}, we deduce from the above that
\begin{equation}
  \label{eq:14}
  \sum_{d\le y}3^{\omega(d)}|R(\Acal,d)|\ll
  \sqrt{yX^*}(Tq)^{3/16+\varepsilon} +\frac{X^*\log^9 (yX^*)}{T}
\end{equation}
for any positive $\varepsilon$.
In particular we can choose $\delta>0$ and take $T=Y(\log q)(\log X^*)^{11}$ 
 together with $y=X/(Y^{1+3/8}q^{\delta+3/8})\ge1$ and obtain
\begin{equation}
  \label{eq:15}
  \sum_{d\le y}3^{\omega(d)}|R(\Acal,d)|\ll_\delta X/(\log X)^2.
\end{equation}
As a conclusion, we can find a $P_2$ in $\Acal$
provided that
\begin{equation}
  \label{eq:16}
  g=\frac{\log X^*}{\log y}\le 2-0.044560
\end{equation}
 i.e.
\begin{equation}
  \label{eq:18}
  X\ge Y^{\frac{11}{8}\frac{3-0.044560}{1-0.044560}}
  q^{(\frac38+\delta)\frac{2-0.044560}{1-0.044560}}.
\end{equation}
We find that
$\frac38\times\frac{2-0.044560}{1-0.044560}=0.7674\dots\le 4/5$ and
$\frac{11}8\times\frac{3-0.044560}{1-0.044560}=4.2532\cdots\le 9/2$.

When $q$ is cube-free, we use the refined bound also recalled in Lemma~\ref{HB}.
\end{proof}

\section{Some additive combinatorics}
\label{Sac}

Our final results involves usage of additive combinatorics, but the
additive combinatorics problem we address can be formulated
independently. We do so, as it leads to a
better understanding of what we do/don't know how to prove.

We thus start with a finite abelian group $G$ written additively.

We also recall for completeness that, when $\Acal$ and $\Bcal$ are two
subsets of some abelian group $G$, the sum $\Acal+\Bcal$ is defined by
\begin{equation*}
  \Acal+\Bcal=\{a+b\ /\  a\in \Acal, b\in \Bcal\}.
\end{equation*}
In particular, the number of representations of a given element is not taken into
account. We use the shortcuts $2\Acal=\Acal+\Acal$ and
$3\Acal=\Acal+\Acal+\Acal$.
The stabilizer, say $H$, of a subset $\Ccal$ of $G$ is the subgroup defined by
\begin{equation*}
  H=\{g\in G\ /\ \forall c\in \Ccal,\ g+c\in \Ccal\}.
\end{equation*}

\begin{lem}\label{maximal}
  A subset $\Ccal$  of a finite abelian group is a union of cosets modulo its stabilizer.
\end{lem}
\begin{proof}
  On denoting by $H$ this stabilizer, it is enough to check that,
  given any element $c$ from $\Ccal$, we indeed have
  $c+H\subset \Ccal$. This follows from the definition of the
  stabilizer. In group theoretical parlance, $\Ccal$ is
  \emph{saturated} modulo~$H$. The lemma follows swiftly.
\end{proof}
\begin{lem}\label{follow}
  Let $\Acal$  be a subset of a finite abelian group
  and let $H$ be the stabilizer of $\Acal + \Acal$. Then $\Acal +
  \Acal+\Acal$ is also a union of cosets modulo~$H$.
\end{lem}
\begin{proof}
  This follows from the equation $\Acal +
  \Acal+\Acal=\Acal +2\Acal=\Acal +(2\Acal+H)$.
\end{proof}

\subsection{Auxiliary lemmas}

\begin{lem}
  \label{SpecialKneser}
  Let $\Acal$  be a subset of a finite abelian group $G$.
  Let $H$ be the stabilizer of $\Acal + \Acal$. Suppose that $\Acal$  meets 
$\lb$ cosets of $H$. Then
$$ |\Acal + \Acal| \geq (2\lb - 1) |H| .$$
\end{lem}
This is \cite[Corollary A.2]{Ramare-Serra-Srivastav*18}, a corollary
of the famous Kneser's Theorem.
While Lemma~\ref{SpecialKneser} is used when the sets we add have a somewhat
small cardinality, the next lemma is tailored for very large sets. 
\begin{lem}
  \label{overcount}
Let $\Acal$ and $\Bcal$ be two subsets of a finite abelian group $G$
satisfying  $|\Acal| + |\Bcal|>|G|$. Then $\Acal + \Bcal = G$.
\end{lem}
This is \cite[Lemma 5.1]{Ramare-Serra-Srivastav*18}.

\begin{lem}
  \label{SpeIndices}
  Let $\Acal$ be a subset of a finite abelian group $G$ such that
  $|\Acal|\ge \eta|G|$, with $\eta>1/3$. Define, for any integer $Y$,
  \begin{equation*}
    \lambda(Y)
    =
    \begin{cases}
      \lceil \eta Y\rceil +1&\text{when $Y\equiv 2[3]$ and $2\le Y\le
        1/(3\eta-1)$,}\\
      \lceil \eta Y\rceil &\text{otherwise.}
    \end{cases}
  \end{equation*}
  For any subgroup $H$ of index $Y$, assume $\Acal$ meets at least
  $\lambda(Y)$ cosets. Then $\Acal+\Acal+\Acal=G$.
\end{lem}
As $|\Acal|/|H|\ge \eta Y$, the subset $\Acal$ always meets at least
$\lceil \eta Y\rceil$ cosets modulo~$H$. 
\begin{proof}
  Let $H$ be the stabilizer of $2\Acal$. By Lemma~\ref{SpecialKneser}
  and~\ref{overcount}, we have $3\Acal=G$ as soon as
  $|\Acal|+|\Acal+\Acal|>|G|$, hence as soon as
  \begin{equation*}
    |\Acal|+(2\lambda-1)|H|>|G|.
  \end{equation*}
  Since $|\Acal|/H\ge \eta Y$, this is certainly true if
  $\eta Y+ 2\lambda-1>Y$, and, since $\lambda\ge \eta Y$, this holds
  when $Y > 1/(3\eta-1)$. This explains the change of definition of
  $\lambda(Y)$ according to whether $Y$  is smaller or larger than
  $1/(3\eta-1)$. When $Y$ is larger, we only use the rather trivial
  value $\lambda(Y)=\lceil \eta Y\rceil$.

  Let us now turn to the non-trivial case $Y\le 1/(3\eta-1)$.
  By Lemma~\ref{maximal}, the set $2\Acal$ is
  a union of cosets modulo~$H$, and by
  Lemma~\ref{SpecialKneser}, of at least
  $2\lambda-1$ of them. This implies that $\Acal+\Acal+\Acal$
  contains at least $3\lambda-2$ cosets modulo~$H$.
  
  Let us write $Y=3y+a$ for $a\in\{0,1,2\}$. We have
  $\lceil \eta Y\rceil>y+\frac{a}{3}$, and thus $\lceil \eta
  Y\rceil\ge y+1 = \frac{Y}{3}+\frac{3-a}{3}$. This implies that
  $3\lambda-2\ge Y+1-a$. When $a\in\{0,1\}$, this is a least $Y$,
  meaning that $\Acal+\Acal+\Acal=G$. When $Y\equiv2[3]$, i.e. when
  $a=2$, and $Y\le1/(3\eta-1)$, our hypothesis gives us the better bound
  $\lambda\ge\lceil \eta Y\rceil+1$. This is enough to complete the proof.
\end{proof}
\subsection{Setting the problem and preliminary conclusions}
\label{prelim}

Let $\Acal$ be a subset of $G$ such that
\begin{equation}
  \label{defeta}
  |\Acal|\ge \eta |G|,\quad \eta>1/3.
\end{equation}
Our problem is to find hypotheses that will lead to the conclusion
that $3\Acal=G$.

Having Lemma~\ref{SpecialKneser} in mind, we consider the stabilizer $H$ of
$2\Acal=\Acal+\Acal$. Here are some facts and further definitions:
\begin{enumerate}
\item  $2\Acal$ and $3\Acal$ are both a union of $H$-cosets, by Lemma~\ref{maximal}
and~\ref{follow}.
\item  On setting $G^*=G/H$ and $\Acal^*=\Acal/H$, it is enough to
  show that $3\Acal^*=G^*$. This is true even if $H$ does not
  stabilize $\Acal$, since we readily check that $3\Acal=3\Acal+H=3(\Acal+H)$. 
\item Notice that $2\Acal^*$ has a trivial stabilizer in $G^*$.
\item We set $Y=|G/H|=|G^*|$ and $\lambda=|\Acal/G|$ so that
  $|2\Acal^*|\ge 2\lambda-1$ by Lemma~\ref{SpecialKneser}.

\end{enumerate}

As a conclusion of Lemma~\ref{SpeIndices} and given $\eta$, only a
finite number of values of $Y$, the index of $H$, are to be
considered, and they are all congruent to~2 modulo~3. Once these
general considerations are set, let us turn to the hypotheses we are
ready to assume. Here are the first two:
\begin{enumerate}
\item[${\ensuremath[C_0]}$] The subset $\Acal$ generates $G$. This is an obvious
  hypothesis which is in fact implied by our other assumptions.
\item[${\ensuremath [C_1]}$] The subset $\Acal$ has a non-empty intersection with
  every coset of any subgroup of index~2.
\end{enumerate}

Here are two 
additional series of hypotheses we are considering:
\begin{enumerate}
\item[${\ensuremath[C_2(Y_0)]}$] The subset $\Acal$ intersects every
  subgroup of index at most $Y_0$ and congruent to~2 mod~3.
\item[${\ensuremath [C_3(Y_0)]}$] Let $K$ be a subgroup of $G$ of
  index at most $Y_0$ and congruent to~2 mod~3. The subset
  $\Acal\bigcup 2\Acal$ has a non-empty intersection with every coset
  modulo~$K$.
\end{enumerate}
Please notice that hypothesis $[C_2(Y_0)]$ does not ask anything
concerning the coset $u+K$ when $u\notin K$, while hypothesis
$[C_3(Y_0)]$ ensures that an element of $\Acal$ \emph{or} of $2\Acal$
belongs to it.

\begin{lem}
  \label{SpeAdditive}
  Let $\Acal^*$ be a subset of a finite abelian group $G^*$ satisfying
  $\Acal^*\cup 2\Acal^*=G^*$ and such that the stabilizer of $2\Acal^*$ is
  $\{0\}$. Then either $3\Acal^*=G^*$ or we have the six conditions\par
  \hfill\begin{minipage}[t]{0.45\linewidth}
    \begin{itemize}
    \item $0\notin \Acal^*$ and $\Acal^*=-\Acal^*$,
    \item The stabilizer of $\Acal^*$ is $\{0\}$,
    \item $\Acal^*\cap 2\Acal^*=\emptyset$,
    \end{itemize}
  \end{minipage}
  \begin{minipage}[t]{0.35\linewidth}
    \begin{itemize}
    \item $|\Acal^*|\le (|G^*|+1)/3$,
    \item $3\Acal^*=G^*\setminus\{0\}$.
    \end{itemize}
  \end{minipage}
  \par\noindent
  In any case, when $|G^*|>2$, we have $4\Acal^*=G^*$.
\end{lem}
Note that the third condition tells us that $\Acal^*$ is a sum-free
set of the finite abelian group $G^*$. Moreover,
the fifth one is a consequence of the
first four. Such sets have been studied, for instance in \cite{Yap*75}
and in \cite{Green-Ruzsa*05}. We ran some numerical experiments when
$G^*$ is the cyclic group with $\ell$ elements say to detect the
existence of the second case of the above lemma. We found that the
situation is quite rigid but not completely so. Here are some
examples:
\begin{itemize}
\item $\Acal^*=\{2,3\}\mod 5$ is our basic example.
\item We explored $\ell\le 29$ with Sage \cite{sage} and here are the
  size of the sets $\Acal^*$ that are possible: \smallskip
  
  \noindent\hspace*{-30pt}\begin{tabular}{c|ccccccccccccccc}
    $\ell$&8&11&17&18&19&20&21&22&23&24&25&26&27&28&29\cr
    \hline
    $|\Acal^*|$&3&4&6&6&6&6,7&6&7&8&7,8&8&7,8,9&8&8,9&8,10     
  \end{tabular}

\item A small cardinality example modulo~71 is given by
  $\Ccal\cup(-\Ccal)$ where $\Ccal=\{1,3,5,17,26,30,32\}\mod 71$.
\end{itemize}

\begin{proof}
  If $0\in \Acal^*$, the assumption $\Acal^*\cup 2\Acal^*=G^*$ readily
  implies that $3\Acal^*=G^*$.  Let us consider
  \begin{equation*}
    \Gamma=\bigl\{g/(g-\Acal^*)\cap 2\Acal^*=\emptyset\bigr\}=G^*\setminus 3\Acal^*
  \end{equation*}
  which we assume to be non-empty. Let $g\in \Gamma$. Since
  $g-\Acal^*$ has no intersection with $2\Acal^*$, it has to be
  included within $\Acal^*$. But since the cardinalities of $\Acal^*$
  and of $g-\Acal^*$ are the same, we have $g-\Acal^*=\Acal^*$.  Let
  $g_1$ and $g_2$ be two elements of $\Gamma$. We have
  $g_1-(g_2-\Acal^*)=(g_1-g_2)+\Acal^*$ on the one hand, while
  $g_1-(g_2-\Acal^*)=g_1-\Acal^*=\Acal^*$ on the other one. This means
  that $g_1-g_2$ stabilizes $\Acal^*$, hence it stabilizes $2\Acal^*$,
  which means that it vanishes. We have thus proved that $\Gamma$ is
  reduced to a single point, say $g_0$. We have
  $\Acal^*\cap 2\Acal^*=g_0-\Acal^*\cap 2\Acal^*=\emptyset$.

  We also
  find that $g_0-2\Acal^*=2\Acal^*$ and that
  $2\Acal^*=(g_0-\Acal^*)+(g_0-\Acal^*)=2g_0-2\Acal^*$. This implies
  that $2\Acal^*=g_0+(g_0-2\Acal^*)=g_0+2\Acal^*$, i.e. that $g_0$
  stabilizes $2\Acal^*$. It is thus equal to~0,
  i.e. $3\Acal^*=G^*\setminus\{0\}$.  Furthermore, Kneser's Theorem
  (Lemma~\ref{SpecialKneser} is enough) tells us that
  $|2\Acal^*|\ge 2|\Acal^*|-1$, and since
  $|\Acal^*|+|2\Acal^*|=|G^*|$, we get $|\Acal^*|\le (|G^*|+1)/3$.
  This upper bound also implies that $2|2\Acal^*|\ge
  2(2|G^*|-1)/3>|G^*|$ when $|G^*|>2$. 
\end{proof}
\subsection{Results}

\begin{thm}
  \label{firstAdditive}
  Let $\Acal\subset G$ be a subset of the finite abelian group $G$ that
  is such that $|\Acal|/|G|>2/5$.
  On assuming $[C_0]$ and $[C_1]$, we have $3\Acal=G$.
\end{thm}
The previous two papers \cite{Ramare-Walker*16} and 
\cite{Ramare-Serra-Srivastav*18} relied on this result.
\begin{proof}
  Let us set $\eta=|\Acal|/|G|$. As $1/(3\eta-1)<1/(3\times\frac25-1)=5$,
  Lemma~\ref{SpeIndices} tells us that only subgroups of index~2 may
  give rise to a difficulty, but this is avoided by
  $[C_0]$ and $[C_1]$.
\end{proof}

\begin{thm}
  \label{secondAdditive}
  Let $\Acal\subset G$ be a subset of the finite abelian group $G$ that
  is such that $|\Acal|/|G|>3/8$.
  On assuming $[C_0]$, $[C_1]$, $[C_2(5)]$ and $[C_3(5)]$, we have $3\Acal=G$.
\end{thm}
This is the main novelty of the present paper on the additive
combinatorics side. 
\begin{proof}
  Let us set $\eta=|\Acal|/|G|$. As $1/(3\eta-1)<1/(3\frac38-1)=8$,
  Lemma~\ref{SpeIndices} tells us that only subgroups of index~2 or~5
  may give a difficulty.   The case~$Y=2$ is ruled out by $[C_0]$ and
  $[C_1]$. In case $Y=5$, by $[C_3(5)]$, Lemma~\ref{SpeAdditive}
  applies. It tells us that $H\bigcap \Acal=\emptyset$, which is
  assumed to be false by $[C_2(5)]$.
\end{proof}
\begin{thm}
  \label{thirdAdditive}
  Let $\Acal\subset G$ be a subset of the finite abelian group $G$ that
  is such that $|\Acal|/|G|=\eta>1/3$. Set $ Y_0=1/(3\eta-1)$.
  On assuming $[C_0]$, $[C_1]$, $[C_2( Y_0)]$ and
  $[C_3( Y_0)]$, we have $3\Acal=G$. 
\end{thm}
\begin{rem}
This theorem shows that in our setting when $G=\Z{q}$ and
$\Acal=\{p\le q^{1+\varepsilon}, (p,q)=1\}$, the Lindel\"of hypothesis in
$q$-aspect for the Dirichlet $L$-functions is enough to obtain the
conclusion, as we prove Assumption $[C_3( Y_0)]$ in
Theorem~\ref{P2}, when $q$ is large enough, and the proof of
Theorem~\ref{T1} easily gets adapted to show that $[C_2( Y_0)]$
also holds, again when $q$ is large enough.
\end{rem}
\begin{proof}
  The proof of Theorem~\ref{secondAdditive} is immediately adapted to
  this case.
\end{proof}
The above results can be sharpened when looking more closely at the
structure of small examples. Here is such a sharpening.
\begin{thm}
  \label{fourthAdditive}
  Let $\Acal\subset G$ be a subset of the finite abelian group $G$ that
  is such that $|\Acal|/|G|=\eta>4/11$.
  On assuming $[C_0]$, $[C_1]$, $[C_2(5)]$,
  $[C_3(8)]$ and that the 2-part of $G$ is isomorphic to
  $(\mathbb{Z}/2\mathbb{Z})^r$ for some $r\ge0$, we have $3\Acal=G$. 
\end{thm}
\begin{proof}
  We adapt the proof of Theorem~\ref{thirdAdditive} and readily
  discover that we only have to worry about the case $Y=8$. Our
  assumption imples that $G^*=G/H$ is isomorphic to $(\mathbb{Z}/2\mathbb{Z})^3$. By
  Lemma~\ref{SpeAdditive}, we may have difficulties when
  $\Acal^*=\Acal/H$ has cardinality~3. But since $\Acal^*$ generates
  $G^*$ that may be viewed as a vector space of dimension~3, $\Acal^*$
  is thus a basis of $G^*$, which would contradict the fact that
  $\Acal^*\bigcup 2\Acal^*=G^*$.
\end{proof}

\section{Three primes. Proof of
  Theorems~\ref{mainthmasymp}, ~\ref{mainthmasympquatro} and~\ref{mainthmasympbis}}
\label{HBT}

By using the techniques of Iwaniec in \cite[Section
2]{Iwaniec*82} (as noted by H. Mikawa in \cite{Mikawa*91}), one can
prove the next lemma.
\begin{lem}[Iwaniec \cite{Iwaniec*82}]
  \label{iwa}
  Let $\varepsilon\in(0,1/200)$ and $\epsilon'>0$ be given.
  For almost all $a$ modulo $q\ge q_0(\varepsilon,\epsilon')$ and
  $q^{\frac{6}{5}+\varepsilon}\le X$, we have
  \begin{equation*}
  \pi(x;q,a)\le
    \frac{2(1+\epsilon')X/\varphi(q)}{\log ({X}/{q^{3/8}})}.
  \end{equation*}
\end{lem}

\begin{proof}[Proof of Theorem~\ref{mainthmasymp}]
  Let $\mc{P}(y)$ be the set of primes below $y$ that do not divide $q$
and let $\mathcal{A}$ be the image of $\mc{P}(y)$ in $G_q = \Z{q}$. We seek to
show that $\mathcal{A} \cdot \mathcal{A} \cdot \mathcal{A} = G_q$.  We select $X=q^{3/2+\epsilon}$. We compute that
  \begin{equation}
    \label{eq:10}
    \frac{\log (X/q^{3/8})}{2\log X}=
    \frac{3}{8}+
    \frac{\epsilon}{12+8\epsilon}.
  \end{equation}
  We use  Lemma~\ref{iwa}. Outside of a set of density $\epsilon''$,
  we have the stated inequality from which we infer  find that
  \begin{equation}
    \label{eq:19}
    \frac{|\Acal|}{\varphi(q)}\ge
    \frac{3}{8(1+\epsilon')}+
    \frac{\epsilon/(1+\epsilon')}{12+8\epsilon}-\epsilon''>3/8
  \end{equation}
  on selecting $\epsilon'$ and $\epsilon''$ properly. To be precise,
  on assuming $\epsilon\le 1$ and $\epsilon'\le 1/2$, we find that
    \begin{align*}
    \frac{3}{8(1+\epsilon')}+
    \frac{\epsilon/(1+\epsilon')}{12+8\epsilon}-\epsilon''
    &\ge
    \frac{3}{8}-\frac{3\epsilon'}{8}+\frac{\epsilon(1-\epsilon')}{20}-\epsilon''
    \\&\ge
    \frac{3}{8}-\frac{3\epsilon'}{8}+\frac{\epsilon}{40}-\epsilon''.
    \end{align*}
    We select
    $\epsilon'=\epsilon''=\epsilon/100$. Theorem~\ref{thirdAdditive}
    applies, the hypotheses being met by 
    Theorem~\ref{T1} with $Y=2$, then with $Y=5$, and by~Theorem~\ref{P2} for index~$Y=5$.
  \end{proof}
  
\begin{proof}[Proof of Theorem~\ref{mainthmasympquatro}]
  Again,
  let $\mc{P}(y)$ be the set of primes below $y$ that do not divide $q$
and let $\mathcal{A}$ be the image of $\mc{P}(y)$ in $G_q = \Z{q}$. We seek to
show that $\mathcal{A} \cdot \mathcal{A} \cdot \mathcal{A} = G_q$.    We select $X=q^{11/8+\epsilon}$. We compute that
  \begin{equation}
    \label{eq:100}
    \frac{\log (X/q^{3/8})}{2\log X}=
    \frac{4}{11}+
    \frac{12\epsilon}{121+88\epsilon}.
  \end{equation}
  We proceed as for the proof of Theorem~\ref{mainthmasymp} with
  Lemma~\ref{iwa}. We find that
  \begin{equation}
    \label{eq:190}
    \frac{|\Acal|}{\varphi(q)}>4/11.
  \end{equation}
  Theorem~\ref{secondAdditive} applies, the hypotheses again being met
  by Theorem~\ref{T1} with $Y=2$, then with $Y=5$ and by
  Theorem~\ref{P2} for index~$Y=5$. The special hypothesis concerning
  the 2-part of $\Acal$ is ensured by the our assumption on~$q$.
\end{proof}

\begin{proof}[Proof of Theorem~\ref{mainthmasympbis}]
We take $X\ge p^{\frac{6}{5}+\epsilon}$ and $p$ large enough in
terms of $\epsilon$. We find that
\begin{equation*}
  \eta=|\mathcal{A}|/\varphi(p)>\frac{11}{32}.
\end{equation*}
Since $1/(3\times\tfrac{11}{32}-1)=32$, the only $Y$ giving trouble are
congruent to 2~modulo~3 and $\le 32$. Theorem~\ref{thirdAdditive} asks
for $[C_3(32)]$, which is granted by Theorem~\ref{P2}, and for
$[C_2(32)]$. 

We check numerically that this may happen only when
$Y\in\{8,11,14,17,20,23,26,29,32\}$.
Our hypothesis on $p$ excludes all these cases, as $Y$ has to be a
divisor of $\varphi(p)=p-1$.
\end{proof}

\section{Four primes. Proof of Theorem~\ref{mainthmasympter}}
Mikawa in \cite{Mikawa*91} proved the next result (be careful: on
page~31, the theorem is stated with a power $(\log x)^{-A}$ but the $-$
sign faded. The correct result is on page~33, line -3).
\begin{lem}[Mikawa \cite{Mikawa*91}]
  \label{mika}
  Let $\epsilon'>0$ and $A>5$ be given.
  For almost all $a$ modulo $q\ge q_0(\epsilon', A)$, we have
  \begin{equation*}
    \pi(x;q,a)\le
    \frac{(2+\epsilon')x/\varphi(q)}{\log ({x^{2/3}}/{q^{1/9}})},
    \quad x^{6/7}\le q\le x/(\log x)^A.
  \end{equation*}
\end{lem}
\begin{proof}[Proof of Theorem~\ref{mainthmasympter}]
As in the proof of Theorem~\ref{mainthmasymp},
  let $\mc{P}(y)$ be the set of primes below $y$ that do not divide $q$
and let $\mathcal{A}$ be the image of $\mc{P}(y)$ in $G_q = \Z{q}$. We seek to
show that $\mathcal{A} \cdot \mathcal{A} \cdot \mathcal{A}\cdot \mathcal{A} = G_q$.
We set $X=q(\log q)^6$. Hence $q\le X/(\log X)^5$ if $q$ is large
  enough. By Lemma~\ref{mika}, we have
  \begin{equation*}
    |\Acal|/\varphi(q)\ge \frac{\log
      ({x^{2/3}}/{q^{1/9}})}{(2+\epsilon')\log x}-\varepsilon
    \ge\frac{5\log
      (X\log X)}{9(2+\epsilon')\log X}-\varepsilon
    \ge \frac{5}{18}-O(\varepsilon)\ge 0.2777
  \end{equation*}
  since $5/18=0.2777\ldots$ and by selecting $\epsilon$ sufficiently small.
Let $H$ be the stabilizer of $\mathcal{A}\cdot \mathcal{A}$, with
index $Y$ and define $|\Acal/H|=n\ge \ceil{0.2777\cdot Y}$.
By using Lemma~\ref{SpecialKneser}, we get
\begin{equation*}
  \frac{2|\mathcal{A} \cdot \mathcal{A}|}{\varphi(q)}
  \ge 2\frac{2\ceil{0.2777\cdot Y}-1}{Y}.
\end{equation*}

The inequality $\ceil{0.2777\cdot Y}\ge 0.2777Y$ shows that this
quantity is $>1$ when $Y\ge 55$.

When $Y$ is smaller, we first notice that $Y=2$ is ruled out by
Theorem~\ref{T1}. Next, we may apply Theorem~\ref{P2}.  On the one
hand, by Lemma~\ref{SpecialKneser}, we have $|2\Acal/H|\ge 2n-1$,
while, by Theorem~\ref{P2}, we also have
$|\Acal/H|+|2\Acal/H|\ge Y$. This implies that $|2\Acal/H|$ is at
least $Y-n$. Thus $|2\Acal/H|$ is at least $\max(2n-1,Y-n)$ which is
greater than $(2Y-1)/3$. As $2(2Y-1)/3>Y$ when $Y\ge3$, we conclude
that $3\Acal=G_q$, as required. 
\end{proof}
\noindent
\thanks{ {\bf Acknowledgements: }The first and second authors have been partly supported by the Indo-French
  Centre for the Promotion of Advanced Research -- CEFIPRA, project No
  5401-1. The first author
    acknowledges the  financial  support  by the Indian National Science Academy
 through a distinguished professorship. The second author was supported by the joint FWF-ANR project Arithrand: FWF: I 4945-N and ANR-20-CE91-0006.}



\end{document}